\documentclass[12pt,a4paper]{article}
\usepackage[dvips]{epsfig}
\usepackage{amsmath}
\usepackage{amsfonts}
\usepackage{amsthm}
\usepackage{amsbsy}
\usepackage{amsgen}
\usepackage{amscd}
\usepackage{amsopn}
\usepackage{amstext}
\usepackage{amsxtra}

\newtheorem{theorem}{Theorem}

\newtheorem{lemma}{Lemma}
\newtheorem{corollary}{Corollary}

\newcommand {\Z} {\mathbb{Z}}
\newcommand {\R} {\mathbb{R}}
\newcommand {\Q} {\mathbb{Q}}

\newcommand {\PP} {\mathbb P}
\begin{document}

\title{Counting Singular Matrices with Primitive Row Vectors}
\author{Igor Wigman\\ School of Mathematical Sciences\\
Tel Aviv University\\ Tel Aviv 69978, Israel}

\maketitle
\begin{abstract}
We solve an asymptotic problem in the geometry of numbers, where we count the number of singular $n\times n$ matrices where row vectors are primitive and of length at most T. Without the constraint of primitivity, the problem was solved by Y. Katznelson. We show that as $T \rightarrow \infty $, the number is asymptotic to 
$ \frac{(n-1)u_n}{\zeta (n) \zeta(n-1)^{n}}T^{n^{2}-n}\log (T)$ for $n \ge 3$. The 3-dimensional case is the most problematic and we need to invoke an equidistribution theorem due to W. Schmidt.
\end{abstract}


AMS classification numbers: 11H06.

\section{Introduction}

\subsection{}
A basic problem in geometry of numbers is counting integer matrices
with certain additional properties. In this paper we will solve a new
counting problem of this kind. Let us consider the set of singular $n
\times n$ matrices with integer entries. We are interested in the
question how many among these matrices have primitive row vectors,
that is each row is \underline{not} a nontrivial multiple of an
integer vector. We count the matrices according to the maximal allowed
Euclidean length of the rows. Without the constraint of primitivity
the problem of counting such matrices was solved by Y. Katznelson
~\cite{Katznelson}. We will find that for $n \ge 3$ a positive
proportion of integer singular matrices have all rows primitive.  

Let $PN_n (T)$ be the counting function of the set $PM_n (T)$ of $n\times n$ singular integer matrices, all of whose rows are primitive and whose Euclidian length is at most $T$. That is, let
\begin{equation*}
PM_{n}(T)= \{M \in M_n (\Z ) :\det (M)=0,\text{primitive rows} \: v_{i}\in \Z^{n}, \, |v_{i}|\leq T\},
\end{equation*}
and set $PN_n (T) = |PM_n (T) |$. In this paper we will determine the
asymptotic behaviour of $PN_n$, as $T \rightarrow \infty$. 

Define a similar counting function $N_n (T)$, where $N_n$ counts $n \times n$ integer matrices M with rows of length $\le T$, that is $N_n (T) = |M_n (T)|$, where
\begin{equation}
\label{eq:M def}
M_{n}(T)= \{M \in M_n (\Z ) :\det (M)=0,\text{rows} \: v_{i},\, |v_{i}|\leq T\}.
\end{equation}

Y. Katznelson \cite{Katznelson} showed that for $n \ge 3$,
\begin{equation*}
N_n (T) = \frac{(n-1)u_n}{\zeta (n)} T^{n^{2}-n}\log (T) + O(T ^ {n ^2 - n} ),
\end{equation*}
where the constant in the $O$-notation depends only on $n$. Recall that $n!!$ denotes the product of integers $\le n$ of the same parity as $n$. The constant $u_n$ is given by:
\begin{equation} 
\label{eq:u_n}
u_n= \; \begin{cases}
\frac{n}{2}(\frac{2(2\pi)^{m-1}}{(2m-1)!!})^n  \cdot \frac{\pi^m}{m!} , \, &n=2m\\
\frac{n}{2}(\frac{\pi^m}{m!})^n \cdot \frac{2(2\pi)^m}{(2m+1)!!} , \, &n=2m+1
\end{cases}.
\end{equation}

Trivially, $PN_n (T) \le N_n (T)$, and thus $PN_n (T) \ll T^{n^2 - n} \log(T)$. Moreover, $PN_n (T) \gg T^{n^2 - n}$, since we can consider, for example, only matrices M with primitive rows $v_i$, which satisfy $v_n = v_1$ and $|v_i| \le T$. A random vector in $\Z ^n$ is primitive with a positive probability, that is, the number of primitive vectors whose length is at most $T$, is $\gg T^n$. The number of such matrices is obviously $\gg (T^n) ^{n-1} = T^ {n^2 - n}$. Combining the observations of this paragraph, we conclude that
$T^{n^2 - n} \ll PN_n (T) \ll T^{n^2 - n} \log(T)$.

\noindent For $n=2$ an elementary argument shows that
\begin{equation}
\label{eq:result n=2}
PN_2(T) = \frac{2\pi}{\zeta(2)} T ^ 2 + O(T).
\end{equation}

\noindent Our main result is:
\renewcommand{\theenumi} {(\roman{enumi})}
\begin{theorem}
\label{thm:main result}
\begin{enumerate}
\item
\label{it:result n >= 4}
For $n \ge  4$ we have 
\begin{equation*}
PN_n (T) = \frac{(n-1)u_n}{\zeta (n) \zeta(n-1)^{n}}T^{n^{2}-n}\log (T) + O(T^{n^2 -n}).
\end{equation*}
\item
\label{it:result n = 3}
For $n=3$ we have 
\begin{equation*}
PN_3 (T) = \frac{2 u_3}{\zeta (3) \zeta(2)^{3}}T^6\log (T) + O(T^6 \log\log (T)).
\end{equation*}
\end{enumerate}
\end{theorem}

\subsection{} 
Another way to treat our problem is to consider it as a counting
problem of rational points with bounded height on a projective
variety, see~\cite{FMT}.

Let $V \subset \underbrace{\PP ^ {n-1} \times \ldots \times \PP  ^ {n-1}}_{n \, times}$ be the projective variety where the determinant vanishes. The {\it height} of a point $X \in \PP  ^ {n-1}(\Q)$ is defined by 
\begin{equation*}
H'(X) = | \tilde {X} |,
\end {equation*}
where $\tilde {X}$ is a primitive integral point in $\Z ^ {n}$ representing X and $|\cdot |$ is the standard Euclidian norm on $\R ^ n$. Now, for $Y = (Y_1,\, \cdots, \, Y_n) \in V$, define the height 
\begin{equation*}
H(Y) = \max_{1 \le i \le n} H'(Y_i).
\end {equation*}
Then $PN_n (T)$ is the number of points $Y\in V$ of height $H(Y) \le T$.

\subsection{} 
We will present now the main idea in the case $n=3$. If M is an integer singular matrix, then all the rows of M lie in a 2-dimensional lattice $\Lambda \subset \Z^3$. Thus we should count triples of vectors lying in a 2-dimensional sublattice of $\Z^3$ and sum over all such lattices.

For a primitive $\lambda\in \Z^3$ we define a 2-dimensional lattice
$L_\lambda = \{ v\in\Z^3:\: v\bot \lambda \} = \lambda^\bot  $. Denote
the subset of all primitive points in $L_\lambda$, $PL_\lambda =
\{\text{primitive } v\in L_\lambda\} $ and $PL_{\lambda} (T) = \{ v\in
PL_{\lambda} \: : \: |v| \le T \}$. Moreover, denote $PA_\lambda = \{
M \in M_3 (\Z ) : \: \text{rows in } PL_\lambda \}$, and $
PA_{\lambda} (T) = \{ M \in PA_{\lambda} \: : \: |\text{rows}| \le T
\}$. Thus $|PA_\lambda(T)| = |PL_\lambda (T)|^3 $. For $\lambda \ne
\lambda '$ the intersection $L_\lambda \cap L_{\lambda'}$ is a set of
integer points on a line, and thus $|PA_\lambda \cap PA_{\lambda'}|
\le 2^3 = 8$. 
It can be shown that the contribution of such intersections is
negligible, and that 
\begin{equation}\label{sum_L}
PN_n (T) \sim \sum\limits_{|\lambda |\ll T^2}{}
\negthickspace ^{'} |PL_\lambda (T)| ^3
\end{equation}
where the last sum is over
primitive $\lambda \in \Z^3$, such that $L_\lambda$ is ``bounded by T"
(see section \ref{sec:n >= 4}). Now 
\begin{equation}\label{Approx P(T)} 
PL_{\lambda} (T) = \frac{ v_2 }{
  \zeta (2) | \lambda |} T^2 + O(\frac {T \log (T)} {|\lambda_1 |} )
\end{equation}
(see section \ref{sec:points count}), and summing the cube of the main
term of $PL_\lambda (T)$ will give the result. 

A complication in dimension 3 is that for some of the lattices
$L_\lambda$ in the sum \eqref{sum_L}, the error term in 
\eqref{Approx  P(T)} is
asymptotically greater than the main term. Such a phenomenon does not
happen for higher dimensions. In order to show that this phenomenon is
\underline{rare} and the contribution of such lattices is negligible,
we will use an equidistribution theorem of Wolfgang Schmidt
~\cite{Schmidt:distribution} (see theorem \ref{thm:cumulating shmidt
  lattices}).  

\subsection{\underline{Contents:}} We will use some known results of counting integer points in $\Z^2$, or more generally, counting points of a sublattice of $\Z^n$, as well as counting {\it primitive} points in such a sublattice. We will give some basic background on lattices in section \ref{sec:background} and some facts concerning counting lattice points will be given in section \ref{sec:points count}. The goal of sections  \ref{sec:n = 3} and \ref{sec:n >= 4} is to prove cases \ref{it:result n = 3} and \ref{it:result n >= 4} of theorem \ref{thm:main result} respectively. 

\subsection{\underline{Acknowledgement:}}
This work was carried out as part of the author's M.Sc. thesis at Tel Aviv 
University, under the supervision of Prof. Zeev Rudnick. 
The author was supported in part by the Israel Science Foundation, founded by 
the Israel Academy of Science and Humanities. 

\section{Background on lattices}
\label{sec:background}
In this section we will give some basic facts which deal with
sublattices of $\Z^n$. For general background see ~\cite{Siegel}. 

\paragraph{Definition:}
Let $\Lambda$ be a lattice. A basis of $\Lambda$, $\{ \lambda_1,\, \lambda_2,\, \ldots \lambda_m \}$, such that the product of the lengths of the vectors in it is minimized is called {\it reduced}. For such a basis we have:
\begin{equation}
\label{eq:reduced basis inequality}
|\lambda_1| \cdot |\lambda_2 | \cdot \ldots \cdot |\lambda_m | \ll \gg \det (\Lambda).
\end{equation}
A basis that satisfies the last inequality has properties similar to a reduced one.

We say that $\Lambda$ is {\it bounded by $T$}, if it has a reduced
basis consisting of vectors of length at most  $T$. 
If a $k$-dimensional lattice has $k$ linearly independent vectors, all
of length at most $T$, it follows that this lattice is
bounded by $cT$, for some constant $c$ that depends only on the
dimension $k$. In that case we will treat it just as if it was bounded by
$T$, since it will affect only some constants in our upper bounds not
affecting the asymptotic behavior.  

Also, for a lattice $\Lambda$, we will denote 
$$
N_{\Lambda}\left( T \right) = |\{v\in\Lambda:\: |v|\leq T \}|
$$
as well as 
$$
P_{\Lambda}\left( T \right) = |\{v\in\Lambda:\: |v|\leq T , 
\quad  v \; \text{primitive}\}|\;.
$$

An $m$-dimensional lattice $\Lambda \subset \Z^n$ is called {\it
  primitive} if there is no $m$-dimensional lattice $\Lambda^*$
properly containing $\Lambda$. In particular, each vector in any basis
of a primitive lattice is  a primitive vector (the converse is
not necessarily true). The orthogonal lattice $\Lambda ^{\bot}$ of
$\Lambda$ consists of all vectors $v \in \Z^n$, such that $v\cdot u =
0$ for all $u\in \Lambda$. It is a primitive integral lattice of
dimension $n-m$.  

If $\Lambda$ is a primitive lattice, then $(\Lambda^{\bot})^{\bot} =
\Lambda$. Also, in this case, it was shown in ~\cite{Schmidt} (chapter
1, formula (4)), that 
\begin{equation}
\label{eq:orthogonal det}
\det(\Lambda) = \det(\Lambda^{\bot})
\end{equation}

\section{Counting lattice points}
\label{sec:points count}
The goal of this section will be to give some expressions for the number of integer points in a lattice, as well as estimations for the error terms of these expressions, which correspond to primitive lattices.

The next lemma is a basic one, which could be found in different variations in the literature, see e.g. ~\cite{Schmidt}, lemma 2. 

\begin{lemma}
\label{lem:estimate vectors of lattice}
Let $ \Lambda \subset  \R^n$ be an $m$-dimensional lattice, and let $\{ \lambda_1,\, \ldots ,\, \lambda_m\}$ be a reduced basis for $\Lambda$, sorted in increasing order of their norms. Denote 
$|\lambda_i|= \mu_i$ for $1 \le i \le m$. Let $\beta \subset \R^m$ be an $m$-dimensional convex body containing $m$ linearly independent vectors of $\Lambda$, then 
\begin{equation*}
|\Lambda\cap \beta| = \frac{vol\left(\beta\right)}{\det (\Lambda)}+O\left(\frac{vol\left(\partial\beta\right)}{\mu_1 \cdot \ldots \cdot  \mu_{m-1}}\right)
\end{equation*}
(the $\mu' s$ in the denominator are all except for the greatest one).
\end{lemma}

Let $v_n$ be the volume of the standard $n$-dimensional unit ball, that is
\begin{equation*} 
v_n= \; \begin{cases}
\;\;\;\frac{\pi^m}{m!} , \, &n=2m\\
\frac{2(2\pi)^m}{(2m+1)!!} , \, &n=2m+1
\end{cases}.
\end{equation*}

\begin{lemma}
\label{lem:primitive vectors in lattice}
Let $\Lambda\subset \Z^n$ be an $(n-1)$-dimensional primitive lattice
which is bounded by T, for $n \ge 3$. Let $\{\lambda_1,\, \lambda_2
,\ldots\, ,\lambda_{n-1}\}, $ be a  reduced basis of $\Lambda$. Then: 
\renewcommand{\theenumi} {(\roman{enumi})}
\begin{enumerate}
\item
\label{it:primitive vectors in lattice n>=4}
For $n \ge  4$ we have 
\begin{equation*}
P_{\Lambda}\left( T\right) = \frac{v_{n-1}}{\zeta (n-1)\det(\Lambda)}T^{n-1}+O(\frac{T^{n-2}}{|\lambda_1|\cdot |\lambda_2|\cdot\ldots\cdot|\lambda_{n-2}|}).
\end{equation*}
\item
\label{it:primitive vectors in lattice n=3}
For $n=3$ we have 
\begin{equation*}
P_{\Lambda}\left( T\right) = \frac{v_2}{\zeta (2)\det(\Lambda)}T^2+O(\frac{T\log{T}}{|\lambda_1|}).
\end{equation*}
\end{enumerate}
\end{lemma}

\begin{proof}
Since $\Lambda$ is primitive, every vector is a [possibly trivial] integer multiple of a primitive vector \underline{in $\Lambda$}, and thus
$N_{\Lambda}(T) = \sum\limits_{k=1}^{\lfloor T\rfloor}P_{\Lambda}(\frac{T}{k})$, hence by Moebius inversion 
\begin{equation*}
P_{\Lambda}(T)=\sum\limits_{k=1}^{\lfloor T \rfloor}\mu(k)\cdot N_{\Lambda}(\frac{T}{k}).
\end{equation*}
Using on the last expression the result of lemma \ref{lem:estimate vectors of lattice}, where $m=n-1$, and $\beta$ is the $(n-1)$-dimensional  ball on the $(n-1)$-dimensional hyper-plane spanned by $\Lambda$, we get 
\begin{equation*}
\begin{split}
P_{\Lambda}(T)&=\sum\limits_{k=1}^{\lfloor T \rfloor}\mu(k)\cdot \bigl( \frac{v_{n-1}}{\det(\Lambda)}(\frac{T}{k})^{n-1}+O(\frac{T^{n-2}}{k^{n-2} \cdot 
|\lambda_1|\cdot |\lambda_2|\cdot\ldots\cdot|\lambda_{n-2}|}) \bigr) \\
&=\frac{v_{n-1} T^{n-1}}{\det{\Lambda}}\sum\limits_{k=1}^{\lfloor T \rfloor} \bigl( \mu(k){k}^{-(n-1)}+\epsilon_k'(T) \bigr)
=\frac{v_{n-1}}{\zeta (n-1)\det(\Lambda)}T^{n-1}+\epsilon(T). 
\end{split}
\end{equation*}
with  error term $\epsilon(T)$ given by  
\begin{equation*}
\begin{split}
\epsilon(T)&=
\sum\limits_{k=1}^{\lfloor T \rfloor} O(\frac{T^{n-2}}{k^{n-2} \cdot 
|\lambda_1|\cdot |\lambda_2|\cdot\ldots\cdot|\lambda_{n-2}|}) \\
&+ \frac{v_{n-1} \cdot T^{n-1}}{\det{\Lambda}}(\sum\limits_{k=1}^{\lfloor T \rfloor}\mu(k){k}^{-(n-1)}-\frac{1}{\zeta(n-1)}).
\end{split}
\end{equation*}
Thus, since $|\mu(n)|\le 1$ for every $n\in N$,  
\begin{equation*}
\begin{split}
|\epsilon(T)|  &\ll \frac{1}{|\lambda_1| \cdot |\lambda_2|\cdot\ldots\cdot |\lambda_{n-2}|}
\biggl( \sum\limits_{k=1}^{\lfloor T \rfloor }
\frac{T^{n-2}}{k^{n-2}} + T^{n-1}\cdot |\sum\limits_{k=\lfloor T \rfloor + 1}^{\infty}\mu(k){k}^{-(n-1)}|  \biggr) \\
&\ll  \frac{1}{|\lambda_1| \cdot |\lambda_2|\cdot\ldots\cdot |\lambda_{n-2}|} 
\biggl( T^{n-2} \sum \limits_{k=1}^{\lfloor T \rfloor}\frac{1}{k^{n-2}}+T^{n-1}\cdot \frac{1}{T^{n-2}} \biggr) \\
\end{split}
\end{equation*}
We also used here the fact that 
$|\lambda_1|\cdot |\lambda_2|\cdot\ldots\cdot|\lambda_{n-2}| \ll
\det(\Lambda).$ 
Now in order to obtain case \ref{it:primitive vectors in lattice n>=4} of the lemma use the convergence of the series $\sum \limits_{k=1}^{\infty}\frac{1}{k^{n-2}}$ for $n\ge 4$. We use $\sum \limits_ {k=1}^{\lfloor T \rfloor }\frac{1}{k} \ll \log (T)$ in order to prove the other case of the lemma.
\end{proof}

\section{The case $n = 3$}
\label{sec:n = 3}

In this section we will prove case \ref{it:result n = 3} of theorem
\ref{thm:main result}. The computation of the main term is also valid
in the case $n\ge 4$, while for $n=3$ we should be more delicate in order to
obtain the appropriate error term. Thus every lemma which is to be
used in section \ref{sec:n >= 4} will be stated for general $n$ in the
current section. 

The main difficulty in this case is that the error term for estimating $P_{\Lambda}$ (see lemma \ref{lem:primitive vectors in lattice}, case \ref{it:primitive vectors in lattice n=3}) could be asymptotically greater than the main term itself. In order to show that such lattices are rare, and thus their contribution to the error term is asymptotically negligible, we will need an equidistribution result of W. Schmidt ~\cite{Schmidt:distribution} for the space of lattices.

\vspace{5mm} 
Suppose $M$ is a singular matrix. This means that there exists a vector $0 \ne \lambda \in \Z^n$, such that all rows of $M$ are orthogonal to $\lambda$. Thus all the rows of $M$ lie in a $(n-1)$-dimensional lattice $\lambda ^\bot$. Since multiplying $\lambda$ by a constant does not affect this property we can assume that $\lambda$ is primitive. Our basic idea is to sum the number of $n$-tuples of primitive vectors with bounded length lying in $\Lambda$, where $\Lambda \subset \Z^n$ runs over all such lattices. We will see that we can limit the sum to a finite number of such lattices.

Let $\lambda \in \Z^n$ be primitive. Recall that $\lambda ^ {\bot}$ denotes the $(n-1)$-dimensional orthogonal dual lattice to $\lambda$, that is the primitive $(n-1)$-dimensional lattice in $\Z^n$ which consists of all vectors in $\Z^n$ which are orthogonal to $\lambda$. Our discussion leads to the following definitions: denote $L_{\lambda} = \lambda ^ {\bot}$, 

\begin{equation*}
\begin{split}
L_{\lambda} (T) &= \{ v\in L_{\lambda}  \: : \: |v| \le T \}, \\
PL_{\lambda} &= \{ \text{primitive } v\in L_{\lambda} \}, \\
PL_{\lambda}(T) &= PL
_{\lambda} \cap L_{\lambda} (T),
\end{split}
\end{equation*}

Given $\lambda$ denote by $A_{\lambda}$ the set of matrices whose rows lie in $L_{\lambda}$:

\begin{equation}
\label{eq:direct sum}
A_{\lambda} = \{ M\in M_n\left( \Z \right) :\: M\cdot \lambda= 0 \} 
= \underbrace{ L_{\lambda} \oplus L_{\lambda} \oplus \ldots \oplus L_{\lambda}}_{n \, times}, 
\end{equation}

Also denote by $PA_{\lambda}$ the subset of matrices in $A_{\lambda}$, whose rows are primitive. Obviously,
\begin{equation*}
PA_{\lambda} = \underbrace{ PL_{\lambda} \oplus PL_{\lambda} \oplus \ldots \oplus PL_{\lambda}}_{n \, times}
\end{equation*}

Denoting by $A_{\lambda} (T)$ and $PA_{\lambda} (T)$ the set of matrices in $A_{\lambda} $ with (primitive) rows in $L_{\lambda}$ of length $\le T$, we clearly have:
\begin{equation*}
A_{\lambda} \left( T\right)= L_ {\lambda} (T) ^n , 
\end{equation*} 
and 
\begin{equation}
\label{eq:prim direct sum}
PA_{\lambda} \left( T\right)= PL_{\lambda} (T) ^ n.
\end{equation} 

The next lemma connects between the terms just defined with our problem. 
\begin{lemma}[\cite{Katznelson}, Lemma 4]
\label{lem:T-bounded lattice exists}
Let $X\in M_n(T)$. Then there is a primitive $\lambda\in \Z^n$, such that $A_{\lambda}$ is bounded by T and $X\in A_{\lambda}$. 
\end{lemma}

\noindent \underline{Remarks:} \newline
${\bullet}$\indent As it were mentioned, ``bounded by $T$" means bounded by $c_n T$, where the constant $c_n > 0$ depends only on $n$. We will see that $c_n$ doesn't affect our computations, so we will ignore it everywhere except the computation of the main term. \newline
${\bullet}$\indent If $A_\lambda$ is bounded by $T$ then  
$|\lambda| ^ n = \det (A_ \lambda ) \ll T ^ {n^2 -n}$ by \eqref{eq:reduced basis inequality} and \eqref{eq:orthogonal det}, and thus $|\lambda|  \ll T ^ {n-1}$. The converse is \underline{not} true, since there \underline{are} primitive vectors $\lambda \in \Z^n$, such that $|\lambda|\le T^{n-1}$, but $A_{\lambda}$ is not bounded by $T$. We will call such vectors ``bad''; they will only have a minor influence on the asymptotics. \newline
${\bullet}$\indent In every case we will deal with reduced basis, the vectors will be ordered in increasing order of their norms, unless specified otherwise.

We may conclude from lemma \ref{lem:T-bounded lattice exists}, that $M_n \left( T \right) = \bigcup \limits _{|\lambda| \ll T^{n-1}}{}  \negthickspace\negthickspace\negthickspace\negthickspace ^{''} A_{\lambda} \left( T \right)$ and also 
\begin{equation}
\label{eq:set as union}
PM_n \left( T \right) = \bigcup _{|\lambda | \ll T^{n-1}}{} \negthickspace\negthickspace\negthickspace\negthickspace ^{''} PA_{\lambda} \left( T \right).
\end{equation}
Here and everywhere in this paper, we use $\bigcup ^{''}$ and $\sum ^{''}$ to denote a union/sum over primitive vectors $\lambda\in \Z^n$, for which $\lambda ^\bot$ is T-bounded. Analogously, $\bigcup ^{'}$ and $\sum ^{'}$  will denote a union/sum over primitive vectors not saying anything about the orthogonal dual.

It is natural to relate the cardinality of the left side of \eqref{eq:set as union} to sum of the cardinalities of the right side:

\begin{equation*}
\begin{split}
PN_n\left( T\right) & =  
\frac{1}{2}\cdot \sum\limits_{|\lambda|\ll T^{n-1}}{} \negthickspace\negthickspace\negthickspace ^{''} |PA_{\lambda} \left( T \right)| + \epsilon'_1 (T) \\
&=\frac{1}{2}\cdot \sum\limits_{|\lambda|\ll T^{n-1}}{} \negthickspace\negthickspace\negthickspace^{'} |PA_{\lambda} \left( T \right)| + \epsilon'_1 (T) + \epsilon '_3 ( T )\\
\end{split}
\end{equation*}
The factor of $\frac{1}{2}$ is due to the fact that
$PA_{\lambda} = PA_{-\lambda}$;   
the terms  $\epsilon'_1 (T)$ and $\epsilon '_3 \left( T \right)$ are the
error terms which implied by the intersections of $PA_\lambda$ for
different $\lambda$ and the contribution of the so called ``bad"
vectors respectively. (A ``bad" vector is a primitive vector 
$\lambda$, with $|\lambda| \ll T^{n-1}$ such that  $\lambda^{\bot}$ 
is not   bounded by T), which do not allow us to get an estimate of the
primitive vectors contained within it. 
However $|\epsilon '_1 \left( T\right) | \leq |\epsilon_1 \left( T
\right) | $ and
$|\epsilon '_3 \left( T \right) | \leq |\epsilon_3 \left( T \right)|$, 
where $\epsilon_1, \, \epsilon _3 $ are the analogous
error terms in the case of the problem solved in ~\cite{Katznelson},
which were shown to be $O(T^{n^2-n})$ (~\cite{Katznelson} pages
130-133). 
Thus:
\begin{equation}
\label{eq:asymptotic equality to sum}
\begin{split}
PN_n ( T ) &= 
\frac{1}{2} \sum\limits_{|\lambda|\ll T^{n-1}}{}
\negthickspace\negthickspace\negthickspace ^{'} |PA_{\lambda} ( T )|
+O(T^{n^2-n})\\ 
&=
\frac{1}{2}\sum\limits_{|\lambda|\ll T^{n-1}}{}
\negthickspace\negthickspace\negthickspace ^{'} |PL_{\lambda} ( T )| ^
n + O(T^{n^2-n})
\end{split} 
\end{equation}

For $n=3$, we would like to demonstrate how we can achieve the bound for
$\epsilon_1 '$, since in our case it is quite simple. 
Indeed, the only matrices in $PA_{\lambda} \cap PA_{\lambda '} $ with
primitive $\lambda \ne \pm \lambda'$ are of the form $ \begin{pmatrix}
  \pm v \\ \pm v \\ \pm v \end{pmatrix} $ with primitive
$v\in\Z^3$. Conversely, given a primitive $v\in\Z^3$, its contribution
to the sum in \eqref{eq:asymptotic equality to sum} is 
$8 PN_{v^\bot} (c_3 T^2)$, where the 8 factor is the number of all possible
signs of the 3 rows. 
Hence:
\begin{equation*}
|\epsilon_1' (T) | \ll  \sum\limits_{|v| \le T} PN_{v ^\bot} (T^2) \le
\sum\limits_{|v| \le T} N_{v ^\bot} (T^2) \ll
\sum\limits_{|v| \le T}  \frac{T^4}{|v|},
\end{equation*}
Now $\sum\limits_{|v| \le T} \frac{1}{|v|} \ll T^2$
where the last equality is due to summation by parts, making use of the fact that $| \{ \text{primitive } v\in\Z^3 \text{ with } |v| \le T \}| \ll T^3$. Thus $|\epsilon_1 ' (T)| \ll T^6$. 

At this point we would like to substitute the result of lemma
\ref{lem:primitive vectors in lattice} into \eqref{eq:asymptotic
  equality to sum}. This is exactly what we are going to do in case $n
\ge 4$ (see section \ref{sec:n >= 4}). However, for $n=3$, the error
term could be asymptotically greater than the main term. In order to
overcome this difficulty, we will first reduce the last sum to
``convenient'' lattices, that is those with not too big determinant
(lemma \ref{lem:big determinant}) and where the norms of vectors in a
reduced basis do not differ too much (lemma \ref{lem:epicentric
  lattices}). Corollary \ref{cor:estimate good region} will show that
for such lattices the error term is negligible relative to the
corresponding main term. 

\paragraph{}
We will adapt the following notations:
\paragraph{\underline{Notations:}} For a $(n-1)$-dimensional lattice
$\Lambda \subset \Z^n$, we will denote the main term of $(P_{\Lambda}
( T ) )^ n$ as well as the corresponding error term: 
\begin{equation}
\label{eq:not terms lat}
c(T,\, \Lambda) = \frac{v_{n-1} ^ n}{( \det (\Lambda ) ) ^ n \zeta(n-1) ^ n} T^{n^2-n}, \;\;\;\; \epsilon( T , \, \Lambda) = (P_{\Lambda} ( T ) )^ n - c(T,\, \Lambda).
\end{equation}
Moreover, for a vector $\lambda \in \Z ^ n$ denote
\begin{equation}
\label{eq:not terms vec}
c(T,\, \lambda) = c(T,\, \lambda ^ {\bot}), \;\;\;\; \epsilon( T , \, \lambda) = \epsilon( T , \, \lambda^ {\bot}).
\end{equation}

\begin{lemma}
\label{lem:big determinant}
For any constant $A>0$, the following estimate holds:
\begin{equation*}
\sum\limits_{\frac{T^2}{(\log T) ^ A} < \det(\Lambda) \le T^2}{} 
P_{\Lambda}(T)^ 3 \ll T^6 \log\log T,
\end{equation*}
where the sum is over primitive T-bounded lattices $\Lambda \subset \Z^3$.
\end{lemma}
\begin{proof}
We will use here the trivial inequality $P_{\Lambda}(T) \le
N_{\Lambda}(T)$. Now, from lemma  \ref{lem:estimate vectors of
  lattice}, $N_{\Lambda}(T) = \frac{\pi T^2}{\det(\Lambda )} +
O(\frac{T}{|\lambda_1|})  \ll \frac{T^2}{\det(\Lambda)}$, where
$\lambda_1=\lambda_1 (\Lambda)$ is the shortest vector in a reduced
basis of $\Lambda$ (that is, the shortest nontrivial vector in
$\Lambda$). The last inequality is due to $\Lambda$ being bounded by
$T$, since it implies
$\frac{T^2}{\det(\Lambda)} \gg \frac{T^2}{|\lambda_1| \cdot |\lambda_2|} = \frac{T}{|\lambda_1|} \cdot \frac{T}{|\lambda_2|} \ge \frac{T}{|\lambda_1|}$.
We will denote by $n(r)$ the number of primitive two-dimensional
lattices $\Lambda\subset \Z^3$ with $\det(\Lambda) = r$ 
and $N(t) = \sum\limits_{r=1}^{\lfloor t \rfloor}n(r)$. Then
\begin{equation}
\label{eq:lattices estimate}
N(t) \ll t^3,
\end{equation}
since such $\Lambda$ are determined by a primitive vector $\pm
\lambda$ orthogonal to it with $|\lambda|=\det(\Lambda)\leq t$; the
number of such vectors is $\ll t^3$. Thus 
\begin{equation*}
\begin{split}
&\sum\limits_{\frac{T^2}{(\log T) ^ A} < \det(\Lambda) \le T^2}{}
  (P_{\Lambda}(T)) ^ 3 \ll T^6 \sum\limits_{\frac{T^2}{(\log T) ^ A}
    < \det(\Lambda) \le T^2}{} \frac{1}{\det(\Lambda) ^ 3} \\ &\le 
T^6\sum\limits_{r=\lceil\frac{T^2}{(\log T) ^ A}\rceil}^{T^2} \frac{n(r)}{r^3} \ll
T^6 \biggl[ \frac{N(t)}{t^3} \biggr| ^{T^2}_{\frac{T^2}{(\log T) ^ A}}  +
\int_{\frac{T^2}{(\log T) ^ A}}^{T^2}\frac{N(t)}{t^4}dt \biggr] \\ 
&\ll T^6 \cdot \log\log T.
\end{split}
\end{equation*}
We used here summation by parts, substituting (\ref{eq:lattices estimate}) to get an estimate for $N(t)$ in order to obtain the last inequality. This concludes the proof of the lemma. It should be noted, that in addition to what was originally stated, we proved here also the following inequality:
\begin{equation}
\label{eq:big determinant N}
\sum\limits_{\frac{T^2}{(\log(T)) ^ A} < det(\Lambda) \le T^2}{} c (T,\, \Lambda ) \ll T^6 \log(\log(T)),
\end{equation}
where $c(T,\, \Lambda)$ is as in \eqref{eq:not terms lat}.
\end{proof}

We will need the following theorem, which is special case of theorem 5 from ~\cite{Schmidt:distribution}. 
\begin{theorem}
\label{thm:cumulating shmidt lattices}
For $a\geq 1$ let $N(a,\, T)$ be the number of lattices $\Lambda\subset \Z^3$ with successive minima, $\{\mu_1,\,\mu_2\}$ which satisfy $\frac{\mu_2}{\mu_1}\ge a$, and $d(\Lambda) \le T$, then:
\begin{equation*}
N(a,\, T) = const \cdot \arcsin(\frac 1{2a})  T^3 + O( a^{-\frac{1}{2}
}\cdot T^{\frac{5}{2}} ) 
\end{equation*}
\end{theorem}

%

We will use Theorem~\ref{thm:cumulating shmidt lattices}  
in order to prove the following lemma:
\begin{lemma}
\label{lem:epicentric lattices}
For any constants $A>0,\,B>1$, the following estimate holds:
\begin{equation*}
\sum\limits_{\stackrel{\det(\Lambda) < \frac{T^2}{(\log T) ^
      A}}{\frac{|\lambda_2|}{|\lambda_1|} > (\log T)^B}}{} 
P_{\Lambda}(T)^3 \ll \frac{T^6}{ (\log T)^{B-1}},
\end{equation*}
where the sum is over primitive T-bounded lattices $\Lambda \subset \Z^3$.
\end{lemma}
\begin{proof}
We will use summation by parts as well as 
Theorem~\ref{thm:cumulating shmidt lattices} 
to bound the sum. In order to do so we will denote 
\begin{equation*}
m_T(r)=| \{  \Lambda \subset \Z^3,\, \text{2-dim. lattice}: \: \det(\Lambda) = r ,\,
|\lambda_2| / |\lambda_1| > (\log(T))^B \} |.
\end{equation*}
so that  $\sum_{r\leq t} m_T(r) = N((\log T)^B,t)$. 

Using the trivial inequality $P_{\Lambda}(T) \le N_{\Lambda}(T) \ll
\frac{T^2}{\det(\Lambda)}$, as in the proof of lemma \ref{lem:big
  determinant}, we have:
\begin{equation*}
\begin{split}
&\sum\limits_{\stackrel{\det(\Lambda) \le \frac{T^2}{(\log T) ^ A}}{\frac{|\lambda_2|}{|\lambda_1|} > (\log T)^B}}{} (P_{\Lambda}(T)) ^ 3  \ll T^6 \cdot \sum\limits_{\stackrel{\det(\Lambda) \le \frac{T^2}{(\log T) ^ A}}{\frac{|\lambda_2|}{|\lambda_1|} > (\log T)^B}}{} \frac{1}{(\det(\Lambda))^3 } \\
&= T^6 \cdot \sum\limits_{r=1}^{\lfloor \frac{T^2}{(\log T) ^
      A}\rfloor }
\frac{m_T(r)}{r^3} \ll
T^6 \cdot \sum\limits_{r=2}^{\lfloor \frac{T^2}{(\log T) ^ A}\rfloor }\frac{m_T(r)}{r^3} \\
&\ll T^6 \cdot \biggl[ \frac{N((\log T)^B,t)}{t^3} \biggr| ^\frac{T^2}{(\log T) ^ A} _ 2 + 
\int_{2}^{\frac{T^2}{(\log T) ^ A}}\frac{N((\log T)^B,t)}{t^4} dt \biggr] 
\end{split}
\end{equation*}
By Theorem~\ref{thm:cumulating shmidt lattices}, this is  
$\ll T^6/ (\log T)^{B-1} $ as required. 

As in the case of lemma \ref{lem:big determinant}, we proved here also:
\begin{equation}
\label{eq:epicentric lattices N}
\sum\limits_{\stackrel{\det(\Lambda) < \frac{T^2}{(\log T) ^ A}}{\frac{|\lambda_2|}{|\lambda_1|} > (\log T)^B}}{} c (T,\, \Lambda ) \ll \frac{T^6}{ (\log T)^{B-1}},
\end{equation}
\end{proof}

\paragraph{}
\begin{lemma}
\label{lem:good error n=3}
Let $\Lambda \subset \Z^3$ be a 2-dimensional lattice, with a reduced basis $\{ \lambda_1,\, \lambda_2 \}$, such that $\det(\Lambda) \le \frac{T^2}{(\log T)^A}$ and $\frac{|\lambda_2|}{|\lambda_1|} \le (\log T)^B$. Then 
\begin{equation}
\frac{T \log T}{|\lambda_1|} \ll \frac{T^2}{\det(\Lambda)} \cdot \frac{1}{(\log T)^{\frac{A-B}{2} - 1}} 
\end{equation}
\end{lemma}
\begin{proof}
\begin{equation*}
|\lambda_2| = \sqrt{\frac{|\lambda_2|}{|\lambda_1|} \cdot (|\lambda_1| |\lambda_2|)} \ll
\sqrt{(\log T)^B \cdot \det(\Lambda)} \ll \frac{T}{(\log T)^{\frac{A-B}{2}}}.
\end{equation*}
Therefore, 
\begin{equation*}
\begin{split}
\frac{T \log T}{|\lambda_1|} 
&= \frac{T \log T}{|\lambda_1| |\lambda_2|} \cdot |\lambda_2| \ll
\frac{T \log T}{\det(\Lambda)} \cdot \frac{T}{(\log T)^{\frac{A-B}{2}}}  \\
&= \frac{T^2}{\det(\Lambda)} \cdot \frac{1}{(\log T)^{\frac{A-B}{2}-1}},
\end{split}
\end{equation*}
which concludes the proof of the lemma.
\end{proof}

\paragraph{}
We will always want to choose the constants A and B, for which $\frac{A-B}{2} - 1 > 0$,
 since in this case the error term of certain counting function will be asymptotically less than the corresponding main term, as we will notice in the following corollary, which follows immediately from the previous lemma. In fact, we would like to choose constants, that will satisfy 
\begin{equation}
\label{eq:constants condition}
\frac{A-B}{2} - 1 \ge 1, \, A > 0 , \, B > 1
\end{equation}
for example, A = 6, B = 2, so this error term will not affect the general error term. 

Using the proof of lemma \ref{lem:primitive orthogonal vectors} below with case \ref{it:primitive vectors in lattice n=3} of lemma \ref{lem:primitive vectors in lattice} and substituting the result of lemma \ref{lem:good error n=3} we obtain:
\begin{corollary}
\label{cor:estimate good region}
Let $\lambda \in \Z^3$ with $|\lambda| \le \frac{T^2}{(\log T)^A}$, such that  $\frac{|\lambda_2|}{|\lambda_1|} \le (\log T)^B$, for constants $A, \, B,$ which satisfy (\ref{eq:constants condition}). Then, under the notations of \eqref{eq:not terms lat},
\begin{equation*}
\epsilon(T, \lambda) \ll c(T, \lambda) \cdot  \frac{1}{(\log T)^{\frac{A-B}{2}-1}} = o(c(T, \lambda)).
\end{equation*}
\end{corollary}

\paragraph{}
We are now ready to finish the proof of  case \ref{it:result n = 3} of
theorem \ref{thm:main result}. 

Recall \eqref{eq:asymptotic equality to sum} and choose the constants $A, \, B$, which satisfy \eqref{eq:constants condition}. We will ignore the difference between $\ll$ and $\le$, which is not significant for bounding the error terms as well as computation of the main term, as we will see a bit later. Thus:
\begin{equation*}
\begin{split}
&\sum\limits_{|\lambda|\leq T^2}{}  \negthickspace\negthickspace ^{'} |PA_{\lambda} \left( T \right)| = 
\sum\limits_{|\lambda|\leq T^2}{} \negthickspace\negthickspace ^{'} (P_{\lambda^{\bot}} (T))^3 \\
&=\sum\limits_{|\lambda|\leq \frac{T^2}{(\log T)^A}}{} \negthickspace\negthickspace\negthickspace\negthickspace\negthickspace ^{'} (P_{\lambda^{\bot}} (T))^3  + 
\sum\limits_{\frac{T^2}{(\log T)^A} < |\lambda| \le T^2}{} \negthickspace\negthickspace\negthickspace\negthickspace\negthickspace\negthickspace\negthickspace\negthickspace ^{'} (P_{\lambda^{\bot}} (T))^3  \\
&= \sum\limits_{\stackrel{|\lambda|\le \frac{T^2}{(\log T)^A} }
{\frac{|\lambda_2|}{|\lambda_1|} \le (\log T)^B}}{}\negthickspace\negthickspace\negthickspace\negthickspace\negthickspace\negthickspace ^{'}   (P_{\lambda^{\bot}} (T))^3  + 
\sum\limits_{\stackrel{|\lambda|\le \frac{T^2}{(\log T)^A} }
{\frac{|\lambda_2|}{|\lambda_1|} > (\log T)^B}}{}\negthickspace\negthickspace\negthickspace\negthickspace\negthickspace\negthickspace ^{'}  (P_{\lambda^{\bot}} (T))^3  + 
O(T^6 \log(\log T)) \\
& = \sum\limits_{\stackrel{|\lambda|\le \frac{T^2}{(\log T)^A} }
{\frac{|\lambda_2|}{|\lambda_1|} \le (\log T)^B}}{}\negthickspace\negthickspace\negthickspace\negthickspace\negthickspace\negthickspace ^{'}   (P_{\lambda^{\bot}} (T))^3 +
O(T^6 \log(\log T)).
\end{split}
\end{equation*}
We used here lemmas \ref{lem:big determinant} and \ref{lem:epicentric lattices} (recall that $B > 1$ because of \eqref{eq:constants condition}).

Substituting the result of corollary \ref{cor:estimate good region} into the last sum we get:
\begin{equation*}
\begin{split}
&\sum\limits_{\stackrel{|\lambda|\le \frac{T^2}{(\log T)^A} }
{\frac{|\lambda_2|}{|\lambda_1|} \le (\log T)^B}}{}\negthickspace\negthickspace\negthickspace\negthickspace\negthickspace\negthickspace ^{'}   (P_{\lambda^{\bot}} (T)) ^3 = 
\sum\limits_{\stackrel{|\lambda|\le \frac{T^2}{(\log T)^A} }
{\frac{|\lambda_2|}{|\lambda_1|} \le (\log T)^B}}{}\negthickspace\negthickspace\negthickspace\negthickspace\negthickspace\negthickspace ^{'}   (c(T, \lambda)) {\rm }
\big( 1 + O(\frac{1}{(\log T)^{\frac{A-B}{2}-1}}) \bigr)  \\ 
&= \bigl( 1 + O( \frac{1}{(\log T)^{\frac{A-B}{2}-1}}) \bigr) \cdot
 \sum\limits_{\stackrel{|\lambda|\le \frac{T^2}{(\log T)^A} }
{\frac{|\lambda_2|}{|\lambda_1|} \le (\log T)^B}}{}\negthickspace\negthickspace\negthickspace\negthickspace\negthickspace\negthickspace ^{'}   (c(T, \lambda))  \\ 
&= \bigl( 1 + O( \frac{1}{(\log T)^{\frac{A-B}{2}-1}}) \bigr) \cdot
\sum\limits_{|\lambda|\le \frac{T^2}{(\log T)^A}} {}\negthickspace\negthickspace\negthickspace\negthickspace ^{'}   (c(T, \lambda))
+ O(T^6)  \\
&= \bigl( 1 + O( \frac{1}{(\log T)^{\frac{A-B}{2}-1}}) \bigr) \cdot
 \sum\limits_{|\lambda|\le T^2} {}\negthickspace ^{'}   (c(T, \lambda))  + O(T^6 \log(\log T)) .
\end{split}
\end{equation*}
We used here corollary \ref{cor:estimate good region} as well as (\ref{eq:big determinant N}) and (\ref{eq:epicentric lattices N}). It should be noted, that usage of lemmas \ref{lem:big determinant} and \ref{lem:epicentric lattices} as ``black boxes'' is not enough in this case, since each summand of the current series is not necessarily less or equal to the corresponding one in these lemmas (because of the error term).  Substituting \eqref{eq:not terms lat} we obtain (writing $\ll$ rather than $\le$ in the domain of summation for consistency with \eqref{eq:asymptotic equality to sum}).

\begin{equation}
\label{eq:sim to sum n=3}
PN_3(T) \sim \sum\limits_{|\lambda|\ll T^2} {}\negthickspace\negthickspace ^{'}   (c(T, \lambda)) = 
\frac{T^6}{\zeta(2) ^ 3} \cdot \sum\limits_{|\lambda|\ll T^2} {}\negthickspace ^{'}  \frac{v_2 ^3}{|\lambda|^3}
\end{equation}

It remains to compute the inner sum. The same computation holds for $n \ge 4$ and will be used in section \ref{sec:n >= 4}, so we will state it for general $n$.

\begin{lemma}
\label{lem:primitive sum estimate}
\begin{equation}
\label{eq:primitive sum estimate}
\sum\limits_{|\lambda|\le M}{} \negthickspace ^{'}  \frac{v_{n-1}^n}{|\lambda|^n} = \frac{2 \cdot u_n}{\zeta(n)}\log (M)+O(1)
\end{equation}
\end{lemma}
\begin{proof}[Proof of lemma \ref{lem:primitive sum estimate}]
Using Moebius inversion on the set of multiples of each primitive vector separately, we obtain:
\begin{equation}
\label{eq:primitive sum}
\sum\limits_{|\lambda|\leq M}{} \negthickspace ^{'} \frac{v_{n-1}^n}{|\lambda|^n}=\sum\limits_{k=1}^{M}\mu(k)k^{-n}\sum\limits_{1\le |u| \le \lfloor M/k \rfloor}\frac{v_{n-1}^n}{|u|^n} 
\end{equation}

Thus, just as in ~\cite{Katznelson}, we get:
\begin{equation}
\label{eq:polar integral}
\sum\limits_{|\lambda|\leq M}{}\negthickspace ^{'}  \frac{v_{n-1}^n}{|\lambda|^n} = \sum\limits_{k=1}^{M}\mu(k)k^{-n} v_{n-1}^n\int\limits_{1\le |x|
\le M/k}{}\frac{dx}{|x|^n}+O(1).
\end{equation}
The reason that the last equality holds is that 
\begin{equation*}
\begin{split}
|\sum\limits_{1\le |u| \le M/k}\frac{v_{n-1}^n}{|u|^n} - \int\limits_{1\le |x|\le M/k}{}\frac{dx}{|x|^n}  \, | &\ll 
\int\limits_{1\le |x|\le M/k}{}\frac{dx}{|x|^{n+1}} \\ 
&\le \int\limits_{|x| \ge 1}{}\frac{dx}{|x|^{n+1}} < \infty,
\end{split}
\end{equation*}
 and the fact that the series $\sum\limits_{k=1}^{\infty} \mu(k) k^{-n}$ converges absolutely.
Computing the last integral in \eqref{eq:polar integral} in polar coordinates we get: 
\begin{equation*}
\begin{split}
\sum\limits_{|\lambda|\leq M}{}\negthickspace ^{'}   \frac{v_{n-1}^n}{|\lambda|^n} &\sim \sum\limits_{k=1}^{M}\mu(k)k^{-n}\cdot 2u_n\cdot ( \log(M)-\log{k}) \\
&=\frac{2\cdot  u_n}{\zeta(n)}\log(M)+O(1),
\end{split}
\end{equation*}
with $u_n = \frac{v_{n-1}^{n}}{2}\int\limits_{S^{n-1}}{}dx$, where the last integral is computed in the usual spherical [polar] coordinates. Thus, using the formula for $v_n$ as well as the fact that $v_n= \frac{1}{n}\int\limits_{S^{n-1}}{}dx$, yields \eqref{eq:u_n}.

Thus, 
\begin{equation}
\sum\limits_{|\lambda|\leq M}{}\negthickspace ^{'}   \frac{v_{n-1}^n}{|\lambda|^n} \sim 
\frac{2\cdot  u_n}{\zeta(n)}\log{M},
\end{equation}
where the error term is O(1), which yields (\ref{eq:primitive sum estimate}) and completes the proof of lemma \ref{lem:primitive sum estimate}.
\end{proof}

Substituting the result of the last lemma in \eqref{eq:sim to sum n=3}
with $M = c_n \cdot T^2$, where $c_n$ is the constant  implied by the
``$\ll $''-notation in \eqref{eq:asymptotic equality to sum}, will
yields case \ref{it:result n = 3} of theorem \ref{thm:main
  result}. The error term $O(1)$ does not bother us, since after
multiplying it by $const\cdot T^{n^2-n}$, while substituting it in
\eqref{eq:asymptotic equality to sum}, we will get an error term of
$O(T^{n^2-n})$, and adding it to other error terms will not increase
an estimate for the general error term of the asymptotics. As
mentioned before, the constant $c_n$ does not affect the computation,
since we substitute it in a logarithm in any case. 

\section{The case $n \ge 4$.}
\label{sec:n >= 4}

We will need the following lemma:

\begin{lemma}
\label{lem:primitive orthogonal vectors}
Let $\lambda \in \Z^n$ be a primitive vector, such that its orthogonal dual, $\lambda^\bot$ is T-bounded. Let $\{\lambda_1,\, \lambda_2 ,\ldots\, ,\lambda_{n-1}\}$ be a reduced basis of $\lambda^\bot$. Then for $n\ge 4$ we have
\begin{equation*}
\begin{split}
|PA_{\lambda}(T)| &=  \frac{v^n_{n-1}}{\zeta ^n (n-1)|\lambda|^n}T^{n^2-n} \\ &+ 
O(\frac{T^{n^2-n-1}}{|\lambda_1|^n\cdot |\lambda_2|^n\cdot\ldots\cdot|\lambda_{n-2}|^n \cdot |\lambda_{n-1}|^{n-1}}).
\end{split}
\end{equation*}
\end{lemma}

\begin{proof}
Due to \eqref{eq:prim direct sum}, we have $|PA_{\lambda} \left( T \right)|=\left( P_{\lambda ^ {\bot}}\left( T \right) \right) ^ n$, and substituting  the case \ref{it:primitive vectors in lattice n>=4} of lemma \ref{lem:primitive vectors in lattice} in the last equality, as well as the fact that $\det(\lambda^{\bot})=|\lambda|$, because of (\ref{eq:orthogonal det}), we get:
$|PA_{\lambda} \left( T \right)| = (\frac{v_{n-1}}{\zeta (n-1)|\lambda|}T^{n-1}+O(\frac{T^{n-2}}{|\lambda_1|\cdot |\lambda_2|\cdot\ldots\cdot|\lambda_{n-2}|}))^n=(a+b)^n$. 
We will use the binomial formula for the last expression. The first summand (that is, $a^n$) is just the main term in the result of the lemma. Now, since the lattice $\lambda^\bot$ is T-bounded, $a$ is asymptotically greater than $b$ (that is $a \gg b$), and thus the only asymptotically significant summand in the binomial is the second one (that is $n\cdot a^{n-1} b$). The coefficient $n$ is constant, and thus, by \eqref{eq:reduced basis inequality}, this is the error term we just stated in the lemma.
\end{proof}

\paragraph{Notations:} In this section we will use the notations in \eqref{eq:not terms vec} as well.

Substituting the result of  Lemma~\ref{lem:primitive orthogonal
  vectors}  into \eqref{eq:asymptotic equality to sum}, we obtain:
\begin{equation}
\label{eq:sim with e2}
\begin{split}
PN_n\left( T\right) &  =  
 \frac{1}{2}\cdot \sum\limits_{|\lambda|\ll T^{n-1}}{} \negthickspace\negthickspace\negthickspace\negthickspace ^{'}
\bigl( \frac{v^n_{n-1}}{\zeta ^n (n-1)|\lambda|^n}T^{n^2-n} +
 \epsilon(T,\lambda) \bigr)  + O(T^{n^2-n})\\
&=\frac{T^{n^2-n}}{2\cdot \zeta^n(n-1)}\cdot \sum\limits_{|\lambda|\ll
 T^{n-1}}{} \negthickspace\negthickspace\negthickspace ^{'}
 \frac{v^n_{n-1}}{|\lambda|^n} + \epsilon_2 ' (T) + O(T^{n^2-n}),
\end{split}
\end{equation}
with 
\begin{equation}
\label{eq:ep2 as sum}
\epsilon_2 ' (T) \ll \sum\limits_{|\lambda|\le T^{n-1}}{} \negthickspace\negthickspace\negthickspace ^{'} | \epsilon(T,\lambda) |.
\end{equation}
We will prove in the end of the section that $\epsilon_2 ' (T) \ll T^ {n^2 - n}$.

\paragraph{}
Using (\ref{eq:primitive sum estimate}) with $M = c_n \cdot
T^{n-1}$,(where $c_n$ is the constant implied by the ``$\ll
$''-notation in \eqref{eq:asymptotic equality to sum}) in
\eqref{eq:sim with e2} will imply 
$$
PN_{n}(T)  = 
\frac{(n-1)u_n}{\zeta (n) \zeta(n-1)^{n}}T^{n^2-n}\log T+O(T^{n^2-n})
$$  
which ends the proof of theorem 1 \ref{it:result n >= 4}.  

Just as in case $n=3$, the error term $O(1)$ does not bother us. Thus,
$PN_{n}(T) = \frac{(n-1)u_n}{\zeta (n) \zeta(n-1)^{n}}T^{n^2-n}\log T+\epsilon '(T)$, where $\epsilon '(T)$ is the error term of this asymptotics. Accumulating all the error terms we confronted with and assuming $\epsilon_2 ' (T) \ll T^{n^2-n}$ (which we will prove immediately),  will imply $|\epsilon '(T)|\ll T^{n^2-n}$.

\paragraph{Error term:} The only error term which is not less or equal to the corresponding error term in \cite{Katznelson} is $\epsilon'_2$. However in this case, we can immediately bound it given the results of the work that was already done by Y. Katznelson.

Under the notations \eqref{eq:not terms vec}, for $n\ge 4$ we have, due to \eqref{eq:ep2 as sum} and lemma \ref{lem:primitive orthogonal vectors}:
\begin{equation*}
\label{eq:error term as sum}
\epsilon'_2(T) \ll \sum\limits_{|\lambda|\leq T^{n-1}}{} \negthickspace\negthickspace\negthickspace ^{'} | \epsilon(T, \lambda)| \ll 
\sum\limits_{|\lambda|\leq T^{n-1}}{} \negthickspace\negthickspace\negthickspace ^{'} \frac{T^{n^2-n-1}}{|\lambda_1|^n\cdot |\lambda_2|^n\cdot\ldots\cdot|\lambda_{n-2}|^n \cdot |\lambda_{n-1}|^{n-1}}.
\end{equation*}
From the definition of $\sum ^{'} $ and $\sum ^{''}$, it is obvious that $\sum ^{''}  \le \sum ^{'} $ as long as only nonnegative numbers are involved. The inequality 
\begin{equation*}
\sum\limits_{|\lambda|\leq T^{n-1}}{} ^{''} \frac{T^{n^2-n-1}}{|\lambda_1|^n\cdot |\lambda_2|^n\cdot\ldots\cdot|\lambda_{n-2}|^n \cdot |\lambda_{n-1}|^{n-1}} \ll
T^{n^2 - n}
\end{equation*}
was showed in ~\cite{Katznelson} (pages 130-133) in the course of proving that $\epsilon_2 (T) \ll T^{n^2 - n}$, where 
$\epsilon_2 (T)$ is the corresponding error term in the case of $N_n (T)$. \qed


\begin{thebibliography}{99}
 
\bibitem[FMT]{FMT} 
 Jens Franke,  Yuri I. Manin and Yuri Tschinkel, 
{\em Rational points of bounded height on Fano varieties.} 
Invent. Math. 95 (1989), no. 2, 421--435. 

\bibitem[K]{Katznelson} Yonathan R. Katznelson.
\emph{Singular Matrices and a Uniform Bound for Congruence Groups of
  $SL_n(\Z)$.} 
Duke Mathematical Journal, 1993, Vol. 69, No. 1, pages 121-136

\bibitem[SCH]{Schmidt} Wolfgang M. Schmidt.
\emph{Asymptotic Formulae for Point Lattices of Bounded Determinant and Subspaces of Bounded Height.}
Duke Mathematical Journal, 1968, No. 35, pages 327-339

\bibitem[SCD]{Schmidt:distribution} Wolfgang M. Schmidt.
\emph{The Distribution of Sublattices of $\Z^m$.}
Monatshefte fur Mathematik, 1998, Vol. 125, pages 37-81

\bibitem[SGL]{Siegel} Carl Ludwig Siegel.
\emph{Lectures on the Geometry of Numbers.}
Springer-Verlag, Berlin, 1988

\end{thebibliography}
\end{document}